\documentclass[a4paper,12pt]{article}
\usepackage{times, url}
\usepackage{amsmath,amssymb,amsthm,amsfonts}
\usepackage{graphicx}
\usepackage{dcolumn}
\usepackage{bm}
\usepackage{color}

\begin{document}
\numberwithin{equation}{section}
\vspace{10mm}

\begin{center}
{\LARGE \bf New series representations for zeta numbers using polylogarithmic identities in combination
 with a polynomial description of Bernoulli numbers }
\vspace{8mm}

{\large \bf J. Braun$^1$, D. Romberger$^2$, H. J. Bentz$^3$}
\vspace{3mm}

$^1$Department Chemie, Ludwig-Maximilians-Universit\"at \\ 
M\"unchen, 81377 M\"unchen, Germany \\
e-mail: \url{juergen.braun@cup.uni-muenchen.de}
\vspace{2mm}

$^2$ Fakult\"at IV, Hochschule Hannover \\ 
Ricklinger Stadtweg 120, 30459  Hannover, Germany \\
e-mail: \url{detlef.romberger@hs-hannover.de}
\vspace{2mm}

$^3$ Institut f\"ur Mathematik und Angewandte Informatik \\ 
Samelsonplatz 1, 31141 Hildesheim, Germany \\
e-mail: \url{bentz@cs.uni-hildesheim.de}

\end{center}
\vspace{10mm}

\begin{abstract}
With this paper we introduce a new series representation of $\zeta(3)$, which is based on the Clausen representation
of odd integer zeta values. Although, relatively fast converging series based on the Clausen representation exist
for $\zeta(3)$, their convergence behavior is very slow compared to BBP-type formulas, and as a consequence they
are not used for explicit numerical computations. The reason is found in the fact that the corresponding Clausen
function can be calculated analytically for a few rational arguments only, where $x=\frac{1}{6}$ is the smallest one.
Using polylogarithmic identities in combination with a polynomial description of the even Bernoulli numbers, the
convergence behavior of the Clausen-type representation has been improved to a level that allows us to challenge
ultimately all BBP-type formulas available for $\zeta(3)$. We present an explicit numerical comparison between one
of the best available BBP formulas and our formalism. Furthermore, we demonstrate by an explicit computation using
the first four terms in our series representation only that $\zeta(3)$ results with an accuracy of $2*10^{-26}$,
where our computation guarantees on each approximation level for an analytical expression for $\zeta(3)$. 
\end{abstract}

\section{INTRODUCTION}
During the last decades BBP-type formulas have been established as the technique of choice for very fast
digit extraction of mathematical constants, as for example, $\pi$, $ln(2)$, $\zeta(3)$ or $\zeta(5)$
\cite{bail97,bail01,broad98,cham03,bor04,bail10,Ade11}. This is because the corresponding algorithms
are simply to implement, where the need of computer memory is very low and no multiple precision
arithmetic software is required \cite{bail10}. Apart from digit extraction interest has grown in BBP-type
formulas in context with statistical randomness of the digit expansions of polylogarithmic constants
\cite{bail97}. 

However, a shortcoming of BBP-type formulas is that a variety of binary degree-1 and degree-2 formulas exist,
but only a few ternary (base 3) or even higher degree BBP-type formulas for polylogarithmic constants are
known. The reason is found in the strong increase of complexity of polylogarithmic functional equations as
a function of the corresponding binary degree. An example is given by Adegoke \cite{kun13} for a polylogarithmic
functional equation of degree 5, where no BBP-type formula for $\zeta(7)$ or higher odd-integer zeta values
has been discovered so far \cite{kun13}.

Concerning the computation of odd-integer valued zeta numbers the so called Clausen representation of zeta
numbers \cite{Hecke44,Hav03,Bentz94} allows for relatively fast digit extraction, which is not restricted
to $\zeta(3)$ and $\zeta(5)$. We have combined this approach recently with a new polynomial representation of
the Bernoulli numbers in connection with Bendersky's L-numbers \cite{Bendersky33}, which appear in context with
the logarithmic Gamma function \cite{braun15}. As a first application approximate calculations of $\zeta$(3),
$\zeta$(5) and $\zeta$(7) in terms this polynomial representation had been presented, where this computational
procedure is applicable to all $\zeta$-values with integer arguments, as well as to related numbers like
Catalan's constant. Compared to digit extraction via corresponding BBP-type formulas the convergence behavior
is not really competitive because the Clausen functions can be calculated analytically for a few rational
arguments only, where $x=\frac{1}{6}$ is the smallest one. In principle one may argue that the speed up in the
convergence should be significant if one would be able to find smaller real-type arguments which also allows 
one for an analytical computation of the corresponding Clausen function. This is indeed possible by the use
of polylogarithmic ladder identities, which exist for $\zeta(3)$ and $\zeta(5)$. In this contribution we present
a first application to the numerical computation of $\zeta(3)$ where we combine polylogarithmic identities
with our polynomial description of Bernoulli numbers to challenge one of the best available BBP-type formulas
typically used for digit extraction of $\zeta(3)$ \cite{bail01}. We demonstrate by an explicit computation
that a fast computation of $\zeta(3)$ is possible, for example with an accuracy of about $10^{-26}$.
Furthermore, we demonstrate that our approach guarantees for an analytical expression of $\zeta(3)$ independently 
from the requested numerical accuracy. At last we present an explicit numerical computation which shows that
our series representation of $\zeta(3)$ converges more than six orders of magnitude faster compared to the
famous BBP-type formula discovered first by Bailey and coworkers \cite{bail07}. 

The paper is organized as follows: in section 2 we remark on the Clausen representation of odd-integer zeta
numbers and present a first computation of $\zeta(3)$ by use of a well known polylogarithmic ladder identity
for $Li_3\left(\frac{1}{2}\right)$. In section 3 we introduce our polynomial representation of the even Bernoulli
numbers and as a consequence for $\zeta(2n)$, $n \in \mathbb N$. This approach is then combined with the
Clausen representation of $\zeta(3)$ to achieve a fast converging series representation, which guarantees
on each approximation level for an analytical expression of $\zeta(3)$. Furthermore, we demonstrate that the 
iterated use of an appropriate polylogarithmic functional equation for $Li_3\left(x\right)$ allows for a
tremendous speed up of the convergence behavior of our series representation. In section 4 we summarize our
results.  

\section{Clausen representation of zeta numbers}
Well-known for a long time is the famous Euler representation \cite{Hecke44,Hav03} of $\zeta(2n)$
with n $\in \mathbb N$: 
\begin{eqnarray}
\zeta (2n) = (-)^{n+1} B_{2n} \frac{(2\pi)^{2n}}{2(2n)!}~.
\end{eqnarray}
For odd integer numbers, as for example, for $n=3$ one finds \cite{Bentz94,Adam98}:

{\bf \large Lemma 2.1} \\
\begin{eqnarray}
Cl_3(x) = \sum_{n=1}^{\infty}\frac{cos(2\pi nx)}{n^3}~=~\zeta(3)&-&3\pi^2 x^2+2\pi^2 x^2 ln(2\pi |x|)
\nonumber \\ &-& 
8\pi^2 \sum_{n=1}^{\infty} \frac{\zeta(2n)}{2n(2n+1)(2n+2)}x^{2n+2}~,
\end{eqnarray}
with x $\in \mathbb R$. A computation of the Clausen function Cl$_3(x)$ for the argument $x=\frac{1}{6}$
results in \cite{Bentz94}:
\begin{eqnarray}
\zeta(3)~=~\frac{\pi^2}{8}-\frac{\pi^2}{12}ln(\frac{\pi}{3})+\frac{\pi^2}{3}\sum_{n=1}^\infty
\frac{\zeta(2n)}{2n(2n+1)(2n+2)}\left(\frac{1}{6}\right)^{2n}~.
\end{eqnarray}
Equation (2.3) converges rather fast, and obviously the convergence could be improved using smaller arguments
for the Clausen function. Unfortunately, this is not for possible x $\in \mathbb Q$, as for smaller rational
arguments as $x=\frac{1}{6}$ partial sums remain in the computation of the Clausen function, which are not
expressible in terms of $\zeta$(3). The way out is the use of polylogarithmic functions, which are widely used
in so called BBP formulas \cite{broad98,kun13,kun10}. For example, it follows for $Li_3\left(\frac{1}{2}\right)$
\cite{broad98}:
\begin{eqnarray}
Li_3\left(\frac{1}{2}\right)=\sum_{n=1}^{\infty}\frac{1}{n^3}\left(\frac{1}{2}\right)^n=\frac{7}{8}\zeta(3)
+\frac{1}{6}(ln(2))^3-\frac{\pi^2}{12}ln(2)~.
\end{eqnarray}
Reformulating $Li_3\left(\frac{1}{2}\right)$ in the following way:

{\bf \large Lemma 2.2} \\
\begin{eqnarray}
Li_3\left(\frac{1}{2}\right)=\sum_{n=1}^{\infty}\frac{1}{n^3}e^{in\theta}~,
\end{eqnarray}
with $\theta$ = i~ln(2) the computation in terms of the corresponding Clausen functions $Cl_3(x)$ and $Sl_3(x)$
results in:

{\bf \large Lemma 3.2} \\
\begin{eqnarray}
\zeta(3)&=&\frac{2\pi^2}{3} ln(2)-6(ln(2))^2 + \frac{2}{3} (ln(2))^3+4(ln(2))^2 ln(ln(2))
\nonumber \\ &+&
16(ln(2))^2~\sum_{n=1}^{\infty} \frac{(-)^{n+1} \zeta(2n)}{2n(2n+1)(2n+2)} \left(\frac{ln(2)}{2\pi}\right)^{2n}~,
\end{eqnarray}
where $Sl_3(x)$ is known analytically. Using the corresponding expression for $Li_3(\frac{1}{2})$ the new argument
$x=\frac{ln(2)}{2\pi}$ is $x\approx \frac{1}{9}$ instead of $x=\frac{1}{6}$, and as a consequence the convergence
is much faster. Furthermore, this series shows up with an alternating sign which provides some benefit in estimating
the convergence properties. For example, $\zeta(3)$ results from the sum of the explicit terms including the first
term (n=1) from infinite series with an error of $\delta \approx 10^{-07}$.

One may notice that the argument $x=\frac{1}{6}$ works for all Clausen functions $Cl_{2n+1}$ as the well known
identity exists \cite{Bentz94}:
\begin{eqnarray}
Cl_{2n+1}(\frac{\pi}{3}) = \frac{1}{2}\left(1-2^{-2n}\right)\left(1-3^{-2n}\right)\zeta(2n+1)~.
\end{eqnarray}

This procedure is also applicable to $\zeta$(2), with \cite{broad98,kun10}: 
\begin{eqnarray}
Li_2\left(\frac{1}{2}\right)=\frac{\pi^2}{12}-\frac{1}{2}\left(ln(2)\right)^2~.
\end{eqnarray}
It follows then for $\zeta$(2):
\begin{eqnarray}
\zeta(2)=2ln(2)(1-ln(ln(2)))-\frac{1}{2}\left(ln(2)\right)^2-4ln(2)~\sum_{n=1}^{\infty} \frac{(-)^{n+1}
\zeta(2n)}{2n(2n+1)} \left(\frac{ln(2)}{2\pi}\right)^{2n}~. \nonumber \\
\end{eqnarray}

Unfortunately, for zeta values with larger integer argument a similar computation seems not possible because for
$Li_n\left(\frac{1}{2}\right)$ with n $>3$ no closed expressions are known \cite{sit87,zhen07,slo07}. For example,
it follows for $Li_4\left(\frac{1}{2}\right)$:
\begin{eqnarray}
Li_4\left(\frac{1}{2}\right)=\frac{15}{16}\zeta(4)-\frac{7}{8}\zeta(3)ln(2)+\frac{1}{4}\zeta(2)(ln(2))^2
+\sum_{n=1}^{\infty}\frac{(-)^n}{(n+1)^3}H_n~,
\end{eqnarray}
where $H_n$ denotes the ordinary finite harmonic series. For the corresponding infinite series no analytical
expression exists. As a consequence the computational scheme introduced here is applicable to a non-trivial
computation of $\zeta(3)$ only, as $\zeta(2)$ is known from Eq.~(2.1) explicitly.

\section{Explicit calculation of $\zeta(3)$ in terms of polylogarithmic identities}
To further improve the convergence in the calculation of $\zeta(3)$ a polynomial representation of the
Bernoulli numbers will be used \cite{braun15}:

{\bf \large Proposition 2.1} \\
\begin{eqnarray}
\zeta(2n)~=~\frac{\zeta(2)^n}{(2n-1)}\sum_{l=1}^n~(-)^{l+1}\left(\begin{array}{c} 
                               n+2-l\\
                               2
                        \end{array}\right)P^{(l)}(n),
\end{eqnarray}
where the P-polynomials are available from the following recursion relation \cite{braun15}:
\begin{eqnarray}
P^{(n-l+1)}(n)~=~6^n\frac{l-1}{2n-l}\sum_{i=l-1}^{n-1}\frac{P^{i-l+2}(i)}{6^i(2n-2i)}~~~~,l>1~.
\end{eqnarray}
with $P^{(1)}(n)=\frac{1}{n}$. As an example, the next three Polynomials result to:
\begin{eqnarray}
P^{(2)}(n)~=~\frac{3}{2*5}
\end{eqnarray}
\begin{eqnarray}
P^{(3)}(n)~=~\frac{3(21n-43)}{2^3*5^2*7}
\end{eqnarray}
\begin{eqnarray}
P^{(4)}(n)=\frac{63n^2-387n+590}{2^4*5^3*7}
\end{eqnarray} 
With this we have:
\begin{eqnarray}
ln~sin(\pi x)&=&ln(\pi x) \nonumber \\ &-& 
\sum_{n=1}^{\infty}\frac{2\zeta(2)^n}{(2n-1)2n} 
\left(\sum_{l=1}^n~(-)^{l+1}
                       \left(\begin{array}{c} 
                               n+2-l\\
                               2
                        \end{array}\right)
P^{(l)}(n) \right)~x^{2n}~.
\end{eqnarray}
For $\zeta(3)$ it follows then:
\begin{eqnarray}
\zeta(3)~&=&~\frac{2\pi^2}{3} ln(2)-6(ln(2))^2 + \frac{2}{3} (ln(2))^3+4(ln(2))^2 ln(ln(2))
\nonumber \\ &+&192\sum_{i=1}^{\infty}c_i
\left(\frac{ln(2)}{2\sqrt6}\right)^{2i}~,
\end{eqnarray}
with
\begin{eqnarray}
c_i~=~\sum_{n=1}^{\infty}\frac{(-)^{n+1}n(n+1)P^{(i)}(n+i-1)}{(2n+2i-3)(2n+2i-2)(2n+2i-1)(2n+2i)}
\left(\frac{ln(2)}{2\sqrt6}\right)^{2n}~.
\end{eqnarray}
Using furthermore the polylogarithmic identities \cite{Ade11}:
\begin{eqnarray}
Li_3\left(\frac{3}{4}\right)+2Li_3\left(\frac{1}{3}\right)+Li_3\left(\frac{1}{4}\right)&=&\frac{19}{6}\zeta(3)
+\frac{1}{3}(ln(3))^3- \frac{4}{3}(ln(2))^3
\nonumber \\ &-& 
\frac{\pi^2}{3}ln(2)+2ln\left(\frac{4}{3}\right)(ln(2))^2
\end{eqnarray}
and
\begin{eqnarray}
Li_3\left(\frac{1}{3}\right)+\frac{1}{4}Li_3\left(\frac{1}{4}\right)+Li_3\left(\frac{2}{3}\right)&=&\frac{15}{8}\zeta(3)+
\frac{1}{6}(ln(2))^3-\frac{\pi^2}{12}ln(2) 
\nonumber \\ &-& 
\frac{1}{6}\left(ln\left(\frac{3}{2}\right)\right)^3+\frac{1}{2}ln(3)\left(ln\left(\frac{3}{2}\right)\right)^2-
\frac{\pi^2}{6}ln\left(\frac{3}{2}\right)~, \nonumber \\
\end{eqnarray}
together with the following functional equation for $Li_3(x)$: \cite{Ade11} 
\begin{eqnarray}
\frac{7}{4}\zeta(3)&=& \frac{1}{4} Li_3 \left ( \left( \frac{1-x}{1+x} \right)^2 \right)
-2Li_3\left( \frac{1-x}{1+x}\right)+2Li_3(1-x)  
\nonumber \\ &+& 
Li_3\left( \frac{1}{1+x} \right)- \frac{1}{2}Li_3(1-x^2)+\frac{\pi^2}{6}ln(1+x)-\frac{1}{3} (ln(1+x))^3~,
\end{eqnarray}
a new identity results with all three arguments of $Li_3(x)$ much closer to 1:
\begin{eqnarray}
6Li_3\left(\frac{2}{3}\right)&+&3Li_3\left(\frac{3}{4}\right)-Li_3\left(\frac{8}{9}\right)
=\frac{91}{12}\zeta(3)-\frac{\pi^2}{2}ln(2)+\frac{7}{3}(ln(2))^3 - \frac{1}{3}(ln(3))^3
\nonumber \\ &-& 
\frac{1}{3}\left(ln\left(\frac{3}{2}\right)\right)^3-2ln\left(\frac{4}{3}\right)(ln(2))^2 
+ ln(3)\left(ln\left(\frac{3}{2}\right)\right)^2
\nonumber \\ &+&
 \frac{2}{3}\left(ln\left(\frac{4}{3}\right)\right)^3~. \nonumber \\
\end{eqnarray}
In a next step we compute $Li_3(x)$ in terms of the corresponding Clausen function by use of the
polynomial representation of the Bernoulli numbers and with the help of Eq.~(3.1). It follows then
for the polylogarithmic function $Li_3(x)$: 
\begin{eqnarray}
Li_3(x)&=&\zeta(3)-\frac{\pi^2}{6}ln\left(\frac{1}{x}\right) + \frac{1}{12}\left(ln\left(\frac{1}{x}\right)\right)^3
+\frac{3}{4}\left(ln\left(\frac{1}{x}\right)\right)^2-\frac{1}{2}\left(ln\left(\frac{1}{x}\right)\right)^2
ln\left(ln\left(\frac{1}{x}\right)\right) 
\nonumber \\ &-&
24\left(ln\left(\frac{1}{x}\right)\right)^2 
\sum^{\infty}_{n=1}\frac{2(-)^{n+1} \sum\limits^n_{l=1}(-)^{l+1}
                       \left(\begin{array}{c} 
                               n+2-l\\
                               2
                        \end{array}\right)P^{(l)}(n)}{(2n-1)2n(2n+1)(2n+2)}
\left(\frac{ln\left(\frac{1}{x}\right)} {2\sqrt6}\right)^{2n+2}~. \nonumber \\
\end{eqnarray}
For $\zeta(3)$ this gives: 
\begin{eqnarray}
\zeta(3) &=& \frac{2\pi^2}{5}ln(3)-\frac{54}{5}\left(ln\left(\frac{3}{2}\right)\right)^2-
\frac{27}{5}\left(ln\left(\frac{4}{3}\right)\right)^2
+\frac{9}{5}\left(ln\left(\frac{9}{8}\right)\right)^2- \frac{6}{5}\left(ln\left(\frac{3}{2}\right)\right)^3 
\nonumber \\ &+&
\left(ln\left(\frac{4}{3}\right)\right)^3+\frac{1}{5}\left(ln\left(\frac{9}{8}\right)\right)^3+ \frac{28}{5}(ln(2))^3-
\frac{4}{5}(ln(3))^3-\frac{4}{5}\left(ln\left(\frac{3}{2}\right)\right)^3 
\nonumber \\ &+&
\frac{36}{5}\left(ln\left(\frac{3}{2}\right)\right)^2ln\left(ln\left(\frac{3}{2}\right)\right)
+\frac{18}{5}\left(ln\left(\frac{4}{3}\right)\right)^2ln\left(ln\left(\frac{4}{3}\right)\right) 
\nonumber \\ &-&
\frac{6}{5}\left(ln\left(\frac{9}{8}\right)\right)^2ln\left(ln\left(\frac{9}{8}\right)\right)
\nonumber \\ &+&
\frac{1728}{5}\sum_{i=1}^{\infty} a_i \left(\frac{ln(\frac{3}{2})}{2\sqrt6}\right)^{2i}+
\frac{864}{5}\sum_{i=1}^{\infty} b_i \left(\frac{ln(\frac{4}{3})}{2\sqrt6}\right)^{2i}-
\frac{288}{5}\sum_{i=1}^{\infty} c_i \left(\frac{ln(\frac{9}{8})}{2\sqrt6}\right)^{2i}~, \nonumber \\
\end{eqnarray}
with 
\begin{eqnarray}
a_i~=~\sum_{n=1}^{\infty}\frac{(-)^{n+1}n(n+1)P^{(i)}(n)}{(2n+2i-3)(2n+2i-2)(2n+2i-1)(2n+2i)}
\left(\frac{ln(\frac{3}{2})}{2\sqrt6}\right)^{2n}~,
\end{eqnarray}
\begin{eqnarray}
b_i~=~\sum_{n=1}^{\infty}\frac{(-)^{n+1}n(n+1)P^{(i)}(n)}{(2n+2i-3)(2n+2i-2)(2n+2i-1)(2n+2i)}
\left(\frac{ln(\frac{4}{3})}{2\sqrt6}\right)^{2n}~,
\end{eqnarray}
and
\begin{eqnarray}
c_i~=~\sum_{n=1}^{\infty}\frac{(-)^{n+1}n(n+1)P^{(i)}(n)}{(2n+2i-3)(2n+2i-2)(2n+2i-1)(2n+2i)}
\left(\frac{ln(\frac{9}{8})}{2\sqrt6}\right)^{2n}~,
\end{eqnarray}
with the arguments $\frac{ln(\frac{3}{2})}{2\sqrt6}$, $\frac{ln(\frac{4}{3})}{2\sqrt6}$ and
$\frac{ln(\frac{9}{8})}{2\sqrt6}$ for the coefficients a$_i$, b$_i$ and c$_i$. Summing up the first
four terms (i=1,2,3,4) from each of the three infinite series together with the explicit terms for
an approximate computation $\zeta(3)$ follows with an error of $\delta \approx 0.3*10^{-17}$. This
is only two orders of magnitude slower in the convergence when compared, for example, to the
famous BBP formula \cite{bail07}:
\begin{eqnarray}
\zeta(3)&=&\frac{1}{672}\sum_{k=0}^{\infty}\left(\frac{1}{4096}\right)^k
\nonumber \\ &&
\hspace*{1.5cm}
\Bigg[\frac{2048}{(24k+1)^3}-\frac{11264}{(24k+2)^3}-\frac{1024}{(24k+3)^3}+\frac{11776}{(24k+4)^3}
\nonumber \\&&
\hspace*{1.5cm}
-\frac{512}{(24k+5)^3}+\frac{4096}{(24k+6)^3}+\frac{256}{(24k+7)^3}+\frac{3456}{(24k+8)^3}
\nonumber \\ && 
\hspace*{1.5cm}
+\frac{128}{(24k+9)^3}-\frac{704}{(24k+10)^3}-\frac{64}{(24k+11)^3}-\frac{128}{(24k+12)^3}
\nonumber \\ && 
\hspace*{1.5cm}
-\frac{32}{(24k+13)^3}-\frac{176}{(24k+14)^3}+\frac{16}{(24k+15)^3}+\frac{216}{(24k+16)^3}
\nonumber \\ && 
\hspace*{1.5cm}
+\frac{8}{(24k+17)^3}+\frac{64}{(24k+18)^3}-\frac{4}{(24k+19)^3}+\frac{46}{(24k+20)^3}
\nonumber \\ && 
\hspace*{1.5cm}
-\frac{2}{(24k+21)^3}-\frac{11}{(24k+22)^3}+\frac{1}{(24k+18)^3} \Bigg]~.
\end{eqnarray}
Within an iterated use of the functional equation (3.10) better and better approximations can be found.
The slowest convergence is found now by $Li_3$($\frac{2}{3}$). As a consequence we rewrite with the
help of (3.10) $Li_3\left(\frac{2}{3}\right)$. It follows first for the polylogarithmic function:
\begin{eqnarray}
Li_3\left(\frac{2}{3}\right) &=& 8Li_3\left(\sqrt{\frac{2}{3}}\right)-8Li_3\left(\frac{2\sqrt2}{\sqrt2+\sqrt3}\right)
-8Li_3\left(\frac{\sqrt2+\sqrt3}{2\sqrt3}\right)
\nonumber \\ &+& 
2Li_3\left(\frac{4\sqrt6}{5+2\sqrt6}\right)+7\zeta(3)-\frac{2\pi^2}{3}ln\left(\frac{2\sqrt2}{\sqrt2+\sqrt3}\right)
\nonumber \\ &+&
\frac{4}{3} \left[ln\left(\frac{2\sqrt2}{\sqrt2+\sqrt3}\right)\right]^3~, \nonumber \\
\end{eqnarray}
and finally $Li_3\left(\frac{2}{3}\right)$ results to: 
\begin{eqnarray}
Li_3\left(\frac{2}{3}\right)&=&\zeta(3)-\frac{2\pi^2}{3}ln\left(\frac{3}{2}\right)
+\frac{1}{12}\left(ln\left(\frac{3}{2}\right)\right)^3+ \frac{3}{2}\left(ln\left(\frac{3}{2}\right)\right)^2
\nonumber \\ &-&
\left(ln\left(\frac{3}{2}\right)\right)^2 ln\left(\frac{1}{2}ln\left(\frac{3}{2}\right)\right)
+\frac{4\pi^2}{3}ln\left(\frac{\sqrt2+\sqrt3}{2\sqrt2}\right)
-\frac{3}{4}\left(ln\left(\frac{\sqrt2+\sqrt3}{2\sqrt2}\right)\right)^3
\nonumber \\ &-&
6\left(ln\left(\frac{\sqrt2+\sqrt3}{2\sqrt2}\right)\right)^2+4\left(ln\left(\frac{\sqrt2+\sqrt3}{2\sqrt3}\right)\right)^2
ln\left(ln\left(\frac{\sqrt2+\sqrt3}{2\sqrt3}\right)\right)
\nonumber \\ &-&
\frac{\pi^2}{3}ln\left(\frac{5+2\sqrt6}{4\sqrt6}\right) +\frac{1}{6}ln\left(\frac{5+2\sqrt6}{4\sqrt6}\right)^3
+\frac{3}{2}\left(ln\left(\frac{5+2\sqrt6}{4\sqrt6}\right)\right)^2
\nonumber \\ &-&
\left(ln\left(\frac{5+2\sqrt6}{4\sqrt6}\right)\right)^2
ln\left(ln\left(\frac{5+2\sqrt6}{4\sqrt6}\right)\right)
+\frac{2\pi^2}{3}ln\left(\frac{2\sqrt3}{\sqrt2+\sqrt3}\right)
\nonumber \\ &+&
\frac{2}{3}\left(ln\left(\frac{2\sqrt3}{\sqrt2+\sqrt3}\right)\right)^3
-6\left(\ln\left(\frac{2\sqrt3}{\sqrt2+\sqrt3}\right)\right)^2
\nonumber \\ &+&
4\left(ln\left(\frac{2\sqrt3}{\sqrt2+\sqrt3}\right)\right)^2
ln\left(ln\left(\frac{2\sqrt3}{\sqrt2+\sqrt3}\right)\right)
\nonumber \\ &+&
4\left(ln\left(\frac{3}{2}\right)\right)^2\sum_{n=1}^{\infty}\frac{(-)^n
\zeta(2n)}{2n(2n+1)(2n+2)}\left(\frac{ln\left(\frac{3}{2}\right)}{4\pi}\right)^{2n}
\nonumber \\ &-&
16\left(ln\left(\frac{\sqrt2+\sqrt3}{2\sqrt2}\right)\right)^2\sum_{n=1}^{\infty}\frac{(-)^n
\zeta(2n)}{2n(2n+1)(2n+2)}\left(\frac{ln\left(\frac{\sqrt2+\sqrt3}{2\sqrt2}\right)} {2\pi}\right)^{2n}
\nonumber \\ &-&
16\left(ln\left(\frac{2\sqrt3}{\sqrt2+\sqrt3}\right)\right)^2\sum_{n=1}^{\infty}\frac{(-)^n
\zeta(2n)}{2n(2n+1)(2n+2)} \left(\frac{ln\left(\frac{2\sqrt3}{\sqrt2+\sqrt3}\right)} {2\pi}\right)^{2n}
\nonumber \\ &+&
4\left(ln\left(\frac{5+2\sqrt6}{4\sqrt6}\right)\right)^2\sum_{n=1}^{\infty}\frac{(-)^n
\zeta(2n)}{2n(2n+1)(2n+2)} \left(\frac{ln\left(\frac{5+2\sqrt6}{4\sqrt6}\right)} {2\pi}\right)^{2n}~, \nonumber \\
\end{eqnarray}
where now four infinite series appear in the computation of $Li_3\left(\frac{2}{3}\right)$. The slowest
convergence is found in the first infinite series with the argument $x=\frac{ln\left(\frac{3}{2}\right)}{4\pi}
\approx \frac{1}{31}$. The other three arguments are much smaller, at least by a factor of two. Inserting now
in each of the four infinite series the polynomial representation of the even Bernoulli numbers (Eq.~(3.1))
the combination of polylogarithmic identities for $Li_3(x)$ with a polynomial description of Bernoulli numbers
has been established, where the polynomial representation guarantees for an additional speed up in the convergence
behavior of all of the four infinite series by more than an order of magnitude.
Finally, at this approximation level $\zeta(3)$ results to:

\begin{eqnarray}
\zeta(3) &=& \frac{48\pi^2}{5}ln\left(\frac{3}{2}\right)+\frac{6\pi^2}{5}ln\left(\frac{4}{3}\right)
-\frac{2\pi^2}{5}ln\left(\frac{9}{8}\right)-\frac{96\pi^2}{5}ln\left(\frac{\sqrt2+\sqrt3}{2\sqrt2}\right)
 \nonumber \\ &-& 
\frac{48\pi^2}{5}ln\left(\frac{2\sqrt3}{\sqrt2+\sqrt3}\right)+\frac{28\pi^2}{5}ln\left(\frac{5+2\sqrt6}{4\sqrt6}\right)
-\frac{108}{5}\left[ln\left(\frac{3}{2}\right)\right]^2
\nonumber \\ &-&
\frac{27}{5}\left[ln\left(\frac{4}{3}\right)\right]^2
+\frac{9}{5}\left[ln\left(\frac{9}{8}\right)\right]^2 +\frac{432}{5}\left[ln\left(\frac{\sqrt2+\sqrt3}{2\sqrt2}\right)\right]^2
\nonumber \\ &+& 
\frac{432}{5}\left[ln\left(\frac{2\sqrt3}{\sqrt2+\sqrt3}\right)\right]^2
-\frac{108}{5}\left[ln\left(\frac{5+2\sqrt6}{4\sqrt6}\right)\right]^2
-\frac{3}{5}\left[ln\left(\frac{4}{3}\right)\right]^3
\nonumber \\ &-&
\frac{6}{5}\left[ln\left(\frac{3}{2}\right)\right]^3
+\frac{1}{5}\left[ln\left(\frac{9}{8}\right)\right]^3+\frac{54}{5}\left[ln\left(\frac{\sqrt2+\sqrt3}{2\sqrt2}\right)\right]^3
\nonumber \\ &-&
\frac{48}{5}\left[ln\left(\frac{2\sqrt3}{\sqrt2+\sqrt3}\right)\right]^3
-\frac{12}{5}\left[ln\left(\frac{5+2\sqrt6}{4\sqrt6}\right)\right]^3
\nonumber \\ &+&
\frac{18}{5}\left[ln\left(\frac{4}{3}\right)\right]^2ln\left[ln\left(\frac{4}{3}\right)\right]
+\frac{72}{5}\left[ln\left(\frac{3}{2}\right)\right]^2ln\left[ln\left(\frac{3}{2}\right)\right]
\nonumber \\ &-&
-\frac{6}{5}\left[ln\left(\frac{9}{8}\right)\right]^2ln\left[ln\left(\frac{9}{8}\right)\right]
-\frac{288}{5}\left[ln\left(\frac{\sqrt2+\sqrt3}{2\sqrt2}\right)\right]^2
ln\left[ln\left(\frac{\sqrt2+\sqrt3}{2\sqrt2}\right)\right]
\nonumber \\ &-&
\frac{288}{5}\left[ln\left(\frac{2\sqrt3}{\sqrt2+\sqrt3}\right)\right]^2
ln\left[ln\left(\frac{2\sqrt3}{\sqrt2+\sqrt3}\right)\right]
\nonumber \\ &+&
\frac{72}{5}\left[ln\left(\frac{5+2\sqrt6}{4\sqrt6}\right)\right]^2 ln\left[ln\left(\frac{5+2\sqrt6}{4\sqrt6}\right)\right]
+\frac{864}{5}\sum_{n=1}^{\infty}c_n^{(1)}\left[\frac{ln\left(\frac{4}{3}\right)}{2\sqrt6}\right]^{2n}
\nonumber \\ &-&
\frac{864}{5}\sum_{n=1}^{\infty}c_n^{(2)}\left[\frac{ln\left(\frac{9}{8}\right)}{2\sqrt6}\right]^{2n}
+\frac{13824}{5}\sum_{n=1}^{\infty}c_n^{(3)}\left[\frac{ln\left(\frac{3}{2}\right)}{4\sqrt6}\right]^{2n}
\nonumber \\ &-& 
\frac{13824}{5}\sum_{n=1}^{\infty}c_n^{(4)}\left[\frac{ln\left(\frac{\sqrt2+\sqrt3}{2\sqrt2}\right)}{2\sqrt6}\right]^{2n}
-\frac{13824}{5}\sum_{n=1}^{\infty}c_n^{(5)}\left[\frac{ln\left(\frac{2\sqrt3}{\sqrt2+\sqrt3}\right)}{2\sqrt6}\right]^{2n}
\nonumber \\ &+&
\frac{3456}{5}\sum_{n=1}^{\infty}c_n^{(6)}\left[\frac{ln\left(\frac{5+2\sqrt6}{4\sqrt6}\right)}{2\sqrt6}\right]^{2n}~. \nonumber \\
\end{eqnarray}
Summing up again the first four terms (n=1,2,3,4) from each of the six infinite series together with all
terms given explicitly $\zeta(3)$ follows now with an error of $\delta \approx 0.37*10^{-21}$. This
is two orders of magnitude faster in the convergence when compared to the BBP formula \cite{bail07}.

It should be mentioned at this stage, that a further advantage of our series representation is that all of these
six types of coefficients can be expressed in terms of elementary functions based on logarithmic expressions. This
allows for a more detailed insight on $\zeta(3)$ as it guarantees on each approximation level an analytical expression
for $\zeta(3)$.

The slowest convergence is now with $Li_3$($\frac{3}{4}$). Therefore, we rewrite the polylogarithmic
function belonging to the argument $\frac{ln\left(\frac{4}{3}\right)}{2\sqrt6}$, again with the help
of the functional equation (3.10). It follows then:
\begin{eqnarray}
Li_3\left(\frac{3}{4}\right) &=& 8Li_3\left(\sqrt{\frac{3}{4}}\right) -8Li_3\left(\frac{2+\sqrt3}{2\sqrt3}\right)
-8Li_3\left(\frac{4}{2+\sqrt3}\right)
\nonumber \\ &+&
2Li_3\left(\frac{7+4\sqrt3}{8\sqrt3}\right)+7\zeta(3)-\frac{2\pi^2}{3}ln\left(\frac{4}{2+\sqrt3}\right)+
\frac{4}{3} \left[ln\left(\frac{4}{2+\sqrt3}\right) \right]^3~, \nonumber \\
\end{eqnarray}
and finally:
\begin{eqnarray}
Li_3\left(\frac{3}{4}\right)&=&\zeta(3)-\frac{2\pi^2}{3}ln\left(\frac{4}{3}\right)
-\frac{1}{12}\left(ln\left(\frac{4}{3}\right)\right)^3
+\frac{3}{2}\left(ln\left(\frac{4}{3}\right)\right)^2-\left(ln\left(\frac{4}{3}\right)\right)^2
ln\left(ln\left(\frac{4}{3}\right)\right)
\nonumber \\ &+&
\frac{4\pi^2}{3}ln\left(\frac{2+\sqrt3}{2\sqrt3}\right)+\frac{2}{3}\left(ln\left(\frac{2+\sqrt3}{2\sqrt3}\right)\right)^3
-6\left(ln\left(\frac{2+\sqrt3}{2\sqrt3}\right)\right)^2
\nonumber \\ &+& 
4\left(ln\left(\frac{2+\sqrt3}{2\sqrt3}\right)\right)^2ln\left(ln\left(\frac{2+\sqrt3}{2\sqrt3}\right)\right)
+2\pi^2ln\left(\frac{4}{2+\sqrt3}\right) + \frac{2}{3}ln\left(\frac{4}{2+\sqrt3}\right)^3
\nonumber \\ &-&
6ln\left(\frac{4}{2+\sqrt3}\right)^2
+4\left(ln\left(\frac{4}{2+\sqrt3}\right)\right)^2ln\left(ln\left(\frac{4}{2+\sqrt3}\right)\right)
-\frac{\pi^2}{3}ln\left(\frac{7+4\sqrt3}{8\sqrt3}\right)
\nonumber \\ &-&
\frac{1}{6}\left(ln\left(\frac{7+4\sqrt3}{8\sqrt3}\right)\right)^3
+\frac{3}{2}\left(ln\left(\frac{7+4\sqrt3}{8\sqrt3}\right)\right)^2
\nonumber \\ &-&
\left(ln\left(\frac{7+4\sqrt3}{8\sqrt3}\right)\right)^2ln\left(ln\left(\frac{7+4\sqrt3}{8\sqrt3}\right)\right)
\nonumber \\ &+&
4\left(ln\left(\frac{4}{3}\right)\right)^2\sum_{n=1}^{\infty}\frac{(-)^n
\zeta(2n)}{2n(2n+1)(2n+2)}\left(\frac{ln\left(\frac{4}{3}\right)}{4\pi}\right)^{2n}
\nonumber \\ &-&
16\left(ln\left(\frac{2+\sqrt3}{2\sqrt3}\right)\right)^2\sum_{n=1}^{\infty}\frac{(-)^n\zeta(2n)}{2n(2n+1)(2n+2)}
\left(\frac{ln\left(\frac{2+\sqrt3}{2\sqrt3}\right)}{2\pi}\right)^{2n}
\nonumber \\ &-&
16\left(ln\left(\frac{4}{2+\sqrt3}\right)\right)^2\sum_{n=1}^{\infty}\frac{(-)^n\zeta(2n)}{2n(2n+1)(2n+2)}
\left(\frac{ln\left(\frac{4}{2+\sqrt3}\right)} {2\pi}\right)^{2n}
\nonumber \\ &+&
4\left(ln\left(\frac{7+4\sqrt3}{8\sqrt3}\right)\right)^2\sum_{n=1}^{\infty}\frac{(-)^n\zeta(2n)}{2n(2n+1)(2n+2)}
\left(\frac{ln\left(\frac{7+4\sqrt3}{8\sqrt3}\right)} {2\pi}\right)^{2n} \nonumber \\
\end{eqnarray}
This procedure can be applied as often as necessary to compute $\zeta(3)$ with a default
accuracy. The only shortcoming is that the number of infinite series increases caused by the
mathematical structure of the functional equation (3.10), where the corresponding arguments appear
as nested roots. As mentioned before, this procedure is not applicable to higher zeta values, as for
example $\zeta(5)$, because no appropriate functional equations exist. In order to finally challenge
the BBP formula \cite{bail07} the polylogarithms $Li_3\left(\sqrt{\frac{2}{3}}\right)$ and
$Li_3\left(\sqrt{\frac{3}{4}}\right)$ appearing with the slowest convergence behavior at this
approximation level will be rewritten with the help of (3.10). It follows:
\begin{eqnarray}
Li_3\left(\sqrt{\frac{2}{3}}\right) &=& 8Li_3\left(\sqrt[4]{\frac{3}{4}}\right)
-8Li_3\left(\frac{\sqrt[4]2+\sqrt[4]3}{2\sqrt[4]2}\right)
-8Li_3\left(\frac{2\sqrt[4]3}{\sqrt[4]2+\sqrt[4]3}\right)
\nonumber \\ &+&
2Li_3\left(\frac{\sqrt2+\sqrt3+2\sqrt[4]6}{4\sqrt[4]6}\right)
+7\zeta(3) -\frac{2\pi^2}{3}ln\left(\frac{2\sqrt[4]3}{\sqrt[4]2+\sqrt[4]3}\right)
\nonumber \\ &+&
\frac{4}{3} \left[ln\left(\frac{2\sqrt[4]3}{\sqrt[4]2+\sqrt[4]3}\right) \right]^3~, \nonumber \\
\end{eqnarray}
and
\begin{eqnarray}
Li_3\left(\sqrt{\frac{3}{4}}\right) &=& 8Li_3\left(\sqrt[4]{\frac{3}{4}}\right)
-8Li_3\left(\frac{\sqrt2+\sqrt[4]3}{2\sqrt[4]3}\right)
-8Li_3\left(\frac{2\sqrt2}{\sqrt2+\sqrt[4]3}\right)
\nonumber \\ &+&
2Li_3\left(\frac{2+\sqrt3+2\sqrt[4]12}{4\sqrt[4]{12}}\right) +7\zeta(3)
-\frac{2\pi^2}{3}ln\left(\frac{2\sqrt2}{\sqrt2+\sqrt[4]3}\right)
\nonumber \\ &+&
\frac{4}{3} \left[ln\left(\frac{2\sqrt2}{\sqrt2+\sqrt[4]3}\right) \right]^3~. \nonumber \\
\end{eqnarray}
To further increase the speed up in the convergence behavior we use the polynom representation of
the Bernoulli numbers where both B$_{2n}$ and B$_{2n-2}$ are involved \cite{braun15}. It follows
then for $Li_3$(x):
\begin{eqnarray}
Li_3(x)&=&\zeta(3)-\frac{\pi^2}{6}ln\left(\frac{1}{x}\right)+
\frac{1}{12}\left(ln\left(\frac{1}{x}\right)\right)^3
+\frac{3}{4}\left(ln\left(\frac{1}{x}\right)\right)^2-\frac{1}{2}\left(ln\left(\frac{1}{x}\right)\right)^2
lnln\left(\frac{1}{x}\right) 
\nonumber \\ && \hspace*{-1.0cm}
-\frac{1}{288}\left(ln\left(\frac{1}{x}\right)\right)^4 + 24\left(ln\left(\frac{1}{x}\right)\right)^2 
\sum^{\infty}_{n=1} \Bigg[
\frac{(-)^{n+1}\sum\limits^{n+1}_{l=1}(-)^{l+1}
                       \left(\begin{array}{c} 
                               n+5-l\\
                               4
                        \end{array}\right)P^{(l)}(n+1)}
                        {(2n-1)2n(2n+1)(2n+2)(2n+3)(2n+4)} \nonumber \\ &-&
\sum^{\infty}_{n=1}\frac{2(-)^{n+1} \sum\limits^n_{l=1}(-)^{l+1}
                       \left(\begin{array}{c} 
                               n+2-l\\
                               2
                        \end{array}\right)P^{(l)}(n)}{(2n-1)2n(2n+1)(2n+2)(2n+3)(2n+4)} \Bigg] 
\left(\frac{ln\left(\frac{1}{x}\right)} {2\sqrt6}\right)^{2n+2} \nonumber \\
\end{eqnarray}
\begin{table}
\begin{center}
\begin{tabular}[t]{||c||c|c|c||}
   \hline
   &&&\\
   $\zeta$(3) & $\zeta$(3)-Zeta series & $\zeta$(3)-(Zeta series+BP+PL) & $\zeta$(3)-BBP formula \cite{bail07} \\
   &&&\\
   \hline
   1st order &&& \\
   n=1 & $\delta$=0.2*10$^{-04}$ &  $\delta$=0.1*10$^{-10}$ & $\delta$=0.7*10$^{-07}$ \\
   \hline
   2nd order &&& \\
   n=2 & $\delta$=0.2*10$^{-06}$ &  $\delta$=0.15*10$^{-15}$ & $\delta$=0.4*10$^{-11}$ \\
   \hline
   3rd order &&& \\
   n=3 & $\delta$=0.3*10$^{-08}$ &  $\delta$=0.2*10$^{-20}$ & $\delta$=0.3*10$^{-15}$ \\
   \hline
   4th order &&& \\
   n=4 & $\delta$=0.4*10$^{-10}$ &  $\delta$=0.2*10$^{-25}$ & $\delta$=0.4*10$^{-19}$ \\
   \hline
\end{tabular}
\caption{Approximate computation of $\zeta$(3) as a function of the summation index n by use of
the Clausen-function representation without and with use of the polynomial representation in
combination with corresponding polylogarithmic identities. The numerical errors are compared
to the BBP-type formula \cite{bail07}}
\end{center}
\end{table}
Calculating all relevant polylogarithms with the formula presented above the additional speed up
in the convergence is more than one order of magnitude. Computing $\zeta$(3) at this approximation
level, again by respecting the first four terms in the corresponding series (n=1,2,3,4), the accuracy
is better than $10^{-25}$. This is more than six orders of magnitude better in the convergence when
compared to \cite{bail07}. The complete numerical comparison with the BBP formula
\cite{bail07} for n=1,2,3 and 4 is presented in Tab. I.

\section{SUMMARY}
In summary, we have presented a unique computational scheme for the explicit calculation of $\zeta$(3) by 
introducing a new series representation of $\zeta(3)$, which is based on the Clausen representation
of odd integer zeta values. By an appropriate combination of polylogarithmic identities with a polynomial
description of the even Bernoulli numbers, we were able to speed up the convergence behavior of the
Clausen-based representation of $\zeta(3)$ to a certain level which is significantly faster than that
of the best BBP-type formulas available for $\zeta(3)$. Furthermore, we have presented a corresponding
numerical comparison between or series representation and one of the best available BBP formulas. Furthermore,
we have demonstrated using the first four terms in our series representation only that $\zeta(3)$ can be
computed with an accuracy of $2*10^{-26}$, where our computation guarantees on each approximation level
for an completely analytical expression for $\zeta(3)$. Finally, we have shown that a computation by use
of the combined polynomial representation of $B_{2n}$ and B$_{2n-2}$ further improves the approximate
calculation of $\zeta(3)$ by more than two orders of magnitude at all approximation levels.


\begin{thebibliography}{99} \expandafter\ifx\csname
natexlab\endcsname\relax\def\natexlab#1{#1}\fi \expandafter\ifx\csname
bibnamefont\endcsname\relax \def\bibnamefont#1{#1}\fi
\expandafter\ifx\csname bibfnamefont\endcsname\relax
\def\bibfnamefont#1{#1}\fi \expandafter\ifx\csname
citenamefont\endcsname\relax \def\citenamefont#1{#1}\fi
\expandafter\ifx\csname url\endcsname\relax \def\url#1{\texttt{#1}}\fi
\expandafter\ifx\csname urlprefix\endcsname\relax\def\urlprefix{URL
}\fi \providecommand{\bibinfo}[2]{#2}
\providecommand{\eprint}[2][]{\url{#2}}

\bibitem{bail97}
Bailey, D. H., Borwein, P. B., Plouffe, S. (1997) On the rapid computation of various polylogarithmic constants,
Mathematics of Computation, 66, 903.

\bibitem{bail01}
Bailey, D. H., Crandall, R. E. (2001) On the random character of fundamental constant expansions, Experimental
Mathematics 10, 175.

\bibitem{broad98}
Broadhurst, D. J. (1998) Polylogarithmic ladders, hypergeometric series and the ten millionth digits of
$\zeta(3)$ and $\zeta(5)$. arXiv:math/9803067v1.

\bibitem{cham03}
Chamberland, M. (2003) Binary BBP-formulae for logarithms and generalized Gaussian-Mersenne primes,
Journal of Integer Sequences 6, 10.

\bibitem{bor04}
Borwein, F., Jonathan, M., Borwein, D., Galway, W., William, F. (2004) Finding and excluding b-ary Machin-type
individual digit formulae, Canad. J. Math. 56, 897.

\bibitem{bail10}
Bailey, D. H. (2010) A compendium of BBP-type formulas for mathematical constants. http://crd.lbl.gov/~dhbailey/
dhbpapers/bbp-formulas.pdf

\bibitem{Ade11}
Adegoke, K. (2011) A novel approach to the discovery of ternary BBP-type formulas for polylogarithm identities,
Notes on Number Theory and Discrete Mathematics, 17, No.1, 4.

\bibitem{kun13}
Adegoke, K. (2013) Formal proofs of degree 5 binary BBP-type formulas, Funct. Approx. Comment. Math. 48, 19.

\bibitem{Hecke44}
Hecke, E. (1944) Herleitung des Euler-Produktes der Zetafunktion und einigerL-Reihen aus ihrer Funktionalgleichung,
Mathematische Annalen, 119, 266.

\bibitem{Hav03}
Havil, J. (2003) Exploring Eulers constant, Princeton University Press.

\bibitem{Bentz94}
Bentz, H. J., Braun, J. (1994)  \"Uber die Werte von $\zeta$(2n+1), Hildesheimer Informatikberichte, 14, 1.

\bibitem{Bendersky33}
L. Bendersky, Sur la fonction Gamma generalisee, Acta Math. {\bf 61} (1933), 263. 

\bibitem{braun15}
Braun, J., Romberger, D., Bentz, H. J. (2015) Fast converging series for zeta numbers in terms of polynomial representations of
Bernoulli numbers, http://www.arXiv.org/abs/math/1503.04636, 1.

\bibitem{bail07}
Bailey, D. H., Borwein, J. M., Calkin, N. J., Girgensohn, R., Luke, D. R., Moll, V. H. (2007)
Experimental Mathematics in Action, Wellesley, MA A K Peters. 

\bibitem{Adam98}
Adamchik, V. S. (1998) Polygamma functions of negative order, J. Comput. Appl. Math., 100, 91. 

\bibitem{kun10}
Adegoke, K. (2010) New binary and ternary digit extraction (BBP-type) formulas for trilogarithm constants,
\textit{New York J. Math.}, 16, 361.

\bibitem{sit87}
Sitaramachandrarao, R (1987) A Formula of S. Ramanujan, Journal of Number Theory, 25, 1. 

\bibitem{zhen07}
De-Yin Zheng, Further summation formulae related to generalized harmonic numbers,
J. Math. Anal. Appl. {\bf 335} (2007), 692. 

\bibitem{slo07}
Zlobin, S. A. (2007) Special values of Generalized Polylogarithms, http://www.arXiv.org/abs/math/0712.1656v1, 1.

\end{thebibliography}
\end{document}